\documentclass[11pt]{amsart}

\theoremstyle{plain}
\newtheorem{thm}{Theorem}[section]
\newtheorem{theorem}[thm]{Theorem}

\newtheorem{lemma}[thm]{Lemma}
\newtheorem{corollary}[thm]{Corollary}
\newtheorem{proposition}[thm]{Proposition}
\theoremstyle{definition}

\newtheorem{notation}[thm]{Notation}

\newtheorem{definition}[thm]{Definition}

\numberwithin{equation}{section}

\newcommand{\p}{\partial}

\newcommand{\sC}{{\mathcal C}}

\newcommand{\sL}{{\mathcal L}}

\newcommand{\sO}{{\mathcal O}}

\newcommand{\sQ}{{\mathcal Q}}

\newcommand{\BO}{{\mathbb O}}

\newcommand{\C}{{\mathbb C}}

\newcommand{\BP}{{\mathbb P}}

\newcommand{\Q}{{\mathbb Q}}

\newcommand{\BS}{{\mathbb S}}

\newcommand{\fg}{{\mathfrak g}}

\newcommand{\fgl}{{\mathfrak g}{\mathfrak l}}

\newcommand{\fso}{{\mathfrak s}{\mathfrak o}}

\newcommand{\aut}{{\mathfrak a}{\mathfrak u}{\mathfrak t}}

\def\Sym{\mathop{\rm Sym}\nolimits}

\def\Hom{\mathop{\rm Hom}\nolimits}

\title[Characterizing subadjoint varieties]{Characterizing subadjoint varieties among Legendrian varieties}

\author{Jun-Muk Hwang}

\thanks{This work was supported by the Institute for Basic Science (IBS-R032-D1).}

\begin{document}
\begin{abstract}
For a symplectic vector space $V$, a projective subvariety $Z \subset \BP V$ is a
Legendrian variety if its affine cone $\widehat{Z} \subset V$ is Lagrangian. In addition to the classical examples of subadjoint varieties associated to simple Lie algebras, many examples of nonsingular Legendrian varieties have been discovered which have positive-dimensional automorphism groups.  We give a characterization of subadjoint varieties among such Legendrian varieties in terms of the isotropy representation. Our proof uses some special features of  the projective third fundamental forms of Legendrian varieties and their relation to the lines on the Legendrian varieties.
\end{abstract}

\maketitle

\noindent {\sc Keywords.} Legendrian variety, projective fundamental form, subadjoint varieties.

\noindent {\sc MSC2020 Classification.} 14M15, 53A20, 53D10

\section{Introduction}\label{s.introduction}

Let $V$ be a complex vector space of dimension $2n+2, n \geq 1$, equipped with a symplectic form $\sigma: \wedge^2 V \to \C$. A projective subvariety $Z \subset \BP V$ is a {\em Legendrian variety} if its affine cone $\widehat{Z} \subset V$ is  Lagrangian  with respect to $\sigma$, namely, it has dimension $n+1$ and its tangent spaces are isotropic with respect to $\sigma$.
Many examples of Legendrian varieties can be produced by Bryant's method (see Section 3 of \cite{Br} or Section 4 of \cite{LM}), which associates to any projective variety   $Z' \subset  \BP^{n+1}$ a Legendrian variety  $ Z \subset \BP^{2n+1}$ birational to $Z'$.  The examples obtained this way are usually singular.
It is  harder to find examples of {\em nonsingular} Legendrian varieties. Classically known nonsingular examples are subadjoint varieties defined as follows.

\begin{definition}\label{d.subadjoint} Let $\fg$ be a simple Lie algebra of Dynkin diagram type different from $A_{\ell}$ or $C_{\ell}.$  Let $X^{\fg} \subset \BP \fg$ be its adjoint variety, namely, the unique closed orbit under the adjoint representation. Then for a base point $o \in X^{\fg}$, the isotropy group at $o$ acts irreducibly on a hyperplane $V^{\fg} \subset T_o X^{\fg}$ and the variety of highest weight vectors of this isotropy action is the subadjoint variety $Z^{\fg} \subset \BP V^{\fg}$.  We can list them explicitly as follows (see Section 1.4.6 of \cite{Hw01} or Table 1 of \cite{Bu06}) according to the type of $\fg$.
\begin{itemize} \item[($\fso_{n+5}$)] The Segre product  $\BP^1 \times \Q^{n-1} \subset \BP^{2n+1}$ of $\BP^1$ and the smooth quadric hypersurface $\Q^{n-1} \subset \BP^{n}$. \item[($G_2$)] The twisted cubic curve $v_3(\BP^1) \subset \BP^3$. \item[($F_4$)]     The Pl\"ucker embedding of the Lagrangian Grassmannian of 3-dimensional isotropic subspaces in  a symplectic vector space of dimension 6.
\item[($E_6$)]   The Pl\"ucker embedding ${\rm Gr}(3; W) \subset \BP \wedge^3 W$  of the Grassmannian  of 3-dimensional subspaces in a  vector space $W$ of dimension 6.
\item[($E_7$)] The spinor embedding $\BS^6 \subset \BP^{31}$ of the spinor variety $\BS^6$ of
 isotropic subspaces of dimension 6 in an orthogonal vector space of dimension 12.
\item[($E_8$)] The 27-dimensional highest weight variety $Z \subset\BP^{55}$ of the basic representation of the exceptional Lie group $E_7$.
\end{itemize} \end{definition}
Conversely, a Legendrian variety $Z \subset \BP V$ which is homogeneous under ${\rm Aut}(Z)$ is isomorphic to one of the  subadjoint varieties (e.g. Theorem 5.11 of \cite{Bu06}).

Bryant's method shows that any compact Riemann surface can be realized as a Legendrian curve in $\BP^3$. Examples of nonhomogeneous nonsingular Legendrian varieties of dimension bigger than 1 have been discovered in \cite{Bu07}, \cite{Bu08}, \cite{Bu09} and \cite{LM}. Especially, Buczy\'nski found several examples  with  $\dim {\rm Aut}(Z) >0$.

A major motivation for studying  nonsingular Legendrian varieties has been its potential application in the classification problem of Fano contact manifolds, all of which are expected to be  homogeneous. One approach to the classification problem of Fano contact manifolds is to use their variety of minimal rational tangents  at a general point, which is known to be a nonsingular Legendrian  variety by Theorem 1.1 of \cite{Ke}. Main Theorem of  \cite{Mok} says that if the variety of minimal rational tangents at a general point of a Fano contact manifold is isomorphic to a subadjoint variety, then the Fano contact manifold is homogenous. From this perspective, it is worthwhile to characterize subadjoint varieties among nonsingular Legendrian varieties in terms of suitable geometric conditions.

In this paper, we give a couple of such characterizations. The first one is in terms of the fundamental forms:

\begin{theorem}\label{t.III}
Let $Z \subset \BP^{2n+1}$ be a nonsingular Legendrian variety of dimension $n \geq 2$.  For a general point $z \in Z$, if the third fundamental form ${\rm III}_{Z,z}$ of $Z$ at $z$ is isomorphic to that of a subadjoint variety, then $Z \subset \BP^{2n+1}$ is isomorphic to the subadjoint variety. \end{theorem}

Each subadjoint variety $Z^{\fg} \subset \BP V^{\fg}$ is a Hermitian symmetric space. When $\fg$  is of exceptional type,  Theorem \ref{t.III} follows from Theorem 1 of \cite{HY}. When $\fg$ is of classical type, we exploit the  geometry lines on the Legendrian variety to prove Theorem \ref{t.III}.

Our next result is in terms of linear automorphisms of $Z \subset \BP V$.

\begin{theorem}\label{t.main}
Let $Z \subset \BP^{2n+1}$ be a nonsingular Legendrian variety of dimension $n \geq 2$, different from a linear subspace $\BP^n \subset \BP^{2n+1}$.   For a general point $z \in Z$, let ${\rm Aut}(Z;z) \subset {\rm PGL}(\C^{n+2})$ be the  group of linear automorphism of $Z$  fixing $z \in Z$ and let $\iota_z: {\rm Aut}(Z;z) \to {\rm GL}(T_z Z)$ be  the isotropy representation on the tangent space at $z.$ Then  $\dim {\rm Ker}(\iota_z) \neq 0$ if and only if $Z$ is a subadjoint variety.  \end{theorem}

Roughly speaking, Theorem \ref{t.main} says that subadjoint varieties are the only nonsingular Legendrian varieties that admit symmetries of higher order.
The proof of Theorem \ref{t.main} uses Theorem \ref{t.III} and a result from \cite{Hw20}
which characterizes the third fundamental forms of the subadjoint varieties. The main issue is how to obtain the required property of the  third fundamental form from the assumption $\dim {\rm Ker}(\iota_z) \neq 0.$ The key result  is  Theorem \ref{t.vecQ}, which gives an explicit description of  the Lie algebra of ${\rm Ker}(\iota_z).$
Theorem \ref{t.vecQ} is of independent interest and seems to be new  even for subadjoint varieties. As
Theorem \ref{t.vecQ} works for a large class of singular Legendrian varieties as well, we expect that it would be useful in the study of the symmetries of  Legendrian varieties.

We remark that there is another characterization of subadjoint varieties   in terms of contact prolongations of ${\rm Aut}(Z)$, proved in \cite{Hw22}. Contact prolongations are completely different from the prolongations of the infinitesimal automorphisms of cubic forms  used below in Section \ref{s.cubic}, and the two methodologies are not directly related.

Let us give a rough outline of the paper. In Section \ref{s.FF}, we review some basic properties of  Legendrian varieties. The proof of Theorem \ref{t.III} is given in Section \ref{s.subadjoint}. In Section \ref{s.cubic}, we recall  a result from \cite{Hw20} on prolongations of cubic forms, with some supplement.   In Section \ref{s.jet}, we prove Theorem \ref{t.vecQ} and obtain Theorem \ref{t.main} as a consequence.

\section{Fundamental forms of Legendrian varieties}\label{s.FF}

Here we collect some results related to fundamental forms of Legendrian varieties.
Some of them are proved in \cite{LM} by the method of moving frames.  Here, we work out the computation in  local coordinates, instead of moving frames.

\begin{notation}\label{n.FF}
Let $V$ be a vector space. Let $Z \subset \BP V$ be a projective subvariety and let $\widehat{Z} \subset V$ be its affine cone.
For a nonsingular point $z\in Z$, we denote by
$${\rm II}_{Z,z} \in \Hom \left(\Sym^2 T_z Z, T_z \BP V/T_z Z \right)$$  the second fundamental form (see page 94 of \cite{FP}) and by $${\rm III}_{Z,z} \in \Hom \left(\Sym^3 T_z Z, T_z \BP V/{\rm Im}({\rm II}_{Z,z}) \right)$$ the  third fundamental form  (see page 126 of \cite{FP}) of $Z$ at $z$.
They are given by the 2-jets (resp. 3-jets) of the restrictions to $Z$ of elements of $H^0(\BP V, \sO(1)) = V^*$ which vanish to the second (resp. third) order at $z$. See Section 2 of \cite{HY} for details.
\end{notation}

We use the following definition of Legendrian varieties, which is equivalent to the one in Section \ref{s.introduction}.

\begin{definition}\label{d.Legendrian}
Let $(V, \sigma)$ be a symplectic vector space of dimension $2n+2, n \geq 1$.
The symplectic form $\sigma$ induces a natural $\sO(2)$-valued 1-form $\theta$ on the projective space $\BP V$ which defines a contact distribution $D$ on $\BP V$
$$ 0 \longrightarrow D \longrightarrow T \BP V \stackrel{\theta}{\longrightarrow} \sO(2) \longrightarrow 0.$$ More precisely, at each point $z\in \BP V$, the hyperplane $D_z \subset T_z \BP V = \Hom(\widehat{z}, V/\widehat{z})$ is defined as $\Hom(\widehat{z}, \widehat{z}^{\perp_{\sigma}}/\widehat{z})$ where $$\widehat{z}^{\perp_{\sigma}} := \{ w \in V \mid \sigma(w, \widehat{z}) =0\}.$$ The distribution $D$ of rank $2n$ is a contact distribution in the sense that ${\rm d} \theta$ induces on $D_z$ a symplectic form for each $z \in \BP V$. A projective subvariety $Z \subset \BP V$ of dimension $n$ is a {\em Legendrian variety} if the tangent space $T_z Z \subset T_z \BP V$ at any nonsingular point $z \in Z$ is contained in $D_z$.
\end{definition}

We  use the following  coordinate system adapted to a Legendrian variety.

\begin{lemma}\label{l.coord}
Let $Z \subset \BP V$ be a  Legendrian variety and let $z \in Z$ be a nonsingular point. Then we can choose an inhomogeneous coordinate system $(x^1, \ldots, x^n, x^{n+1}, \ldots, x^{2n}, x^{2n+1})$ on an affine open subset in $\BP V$ such that \begin{itemize} \item[(1)] $z = (x^1= \cdots = x^{2n+1} =0)$; \item[(2)]  the $n$-dimensional linear space $(x^{n+1} = \cdots = x^{2n+1} =0)$ is tangent to $Z$ at $z$; \item[(3)]
 the contact form $\theta$ is proportional to $$ \sum_{k=1}^n (x^{n+k} {\rm d} x^k - x^k {\rm d} x^{n+k}) - {\rm d} x^{2n+1}. $$ \end{itemize} Furthermore, in a Euclidean neighborhood $U$ of $z$ in $\BP V$, choose coordinates $(y^1, \ldots, y^n)$ on $U \cap Z$ and holomorphic functions $F^1, \ldots, F^n, E$ on $U\cap Z$ satisfying $$F^1(0) = \cdots = F^n(0) = E(0) = 0$$  such that the submanifold
 $Z\cap U$ of $U$ is given by the equations \begin{eqnarray*}  x^k & = & y^k \mbox{ for } 1 \leq k \leq n , \\ x^{n+k} & = & F^k(y^1, \ldots, y^n) \mbox{ for } 1 \leq k \leq n    \mbox{ and } \\ x^{2n+1} & = & E(y^1, \ldots, y^n). \end{eqnarray*}  Then, shrinking $U$ if necessary, we can find a holomorphic function $F$ on $Z \cap U$ such that
 \begin{equation}\label{e.i} F^k = \frac{\partial F}{\partial y^k} \mbox{  for each } 1 \leq k \leq n, \end{equation} \begin{equation}\label{e.ii} \frac{\p E(0)}{\p y^i} = 0 = \frac{ \partial^2 F(0)}{\partial y^i \partial y^k}  \mbox{ for each } 1 \leq i, k \leq n, \end{equation} \begin{equation}\label{e.iii}  - \frac{ \partial^N E}{\partial y^{i_1} \cdots \partial y^{i_N}} = (N-2) \frac{ \partial^N F}{\partial y^{i_1} \cdots \partial y^{i_N}} + \sum_{k=1}^n y^k \frac{ \partial^{N+1} F}{\partial y^k \partial y^{i_1} \cdots \partial y^{i_N}} \end{equation} for each integer $N \geq 1,$ and \begin{equation}\label{e.iv}
 \frac{\p^2 E(0)}{\p y^i \p y^k}  = 0 \mbox{ for each } 1 \leq i, k \leq n. \end{equation}
\end{lemma}

\begin{proof}
The existence of an inhomogeneous coordinate system $$(x^1, \ldots, x^n, x^{n+1}, \ldots, x^{2n}, x^{2n+1})$$ satisfying (1), (2) and (3) is well-known: they arise  from a choice of symplectic basis of $V$ with respect to $\sigma$. Let us verify (\ref{e.i})--(\ref{e.iv}).

Since $Z$ is Legendrian, the pull-back of $\theta$ to $Z$ must vanish:
  \begin{eqnarray}\label{e.E} \sum_{k=1}^n (  F^k {\rm d} y^k - y^k {\rm d} F^k) - {\rm d} E &=& 0. \end{eqnarray} Taking derivative of (\ref{e.E}), we have $ \sum_{k=1}^n {\rm d} F^k \wedge {\rm d} y^k = 0.$ By Poincar\'e lemma,  there exists a holomorphic function $F$ in a small neighborhood of $z$ in $Z$
 such that $F^k = \frac{\partial F}{\partial y^k}$ for each $1 \leq k \leq n$,  verifying (\ref{e.i}).
  Since ${\rm d} E$ and ${\rm d} F^k$ vanish at $z$ by (2), we have (\ref{e.ii}).
    Putting (\ref{e.i}) into  (\ref{e.E}), we have
 \begin{eqnarray*} {\rm d} E  &=& \sum_{k=1}^n (F^k {\rm d} y^k - y^k {\rm d} F^k) \\
 &=& \sum_{k=1}^n \left(\frac{\partial F}{\partial y^k} {\rm d} y^k - y^k \sum_{i=1}^n \frac{\partial^2 F}{\partial y^i \partial y^k} {\rm d} y^i \right) \\ &=& \sum_{k=1}^n  \left( \frac{\partial F}{\partial y^k} - \sum_{i=1}^n y^i \frac{\partial^2 F}{\partial y^i \partial y^k} \right) {\rm d} y^k, \end{eqnarray*} which gives (\ref{e.iii}) for $N=1$.
 Taking higher derivatives   successively, we obtain (\ref{e.iii}) for all $N \geq 2$.
 Finally, (\ref{e.iv}) follows from (\ref{e.iii}).  \end{proof}

\begin{proposition}\label{p.II}
Let $Z \subset \BP V$ be a Legendrian variety and let $z \in Z$ be a nonsingular point. Then \begin{itemize} \item[(a)] ${\rm Im}({\rm II}_{Z,z}) \subset D_z/T_z Z $ and \item[(b)]  ${\rm Im}({\rm II}_{Z,z}) = D_z/T_z Z$ if and only if the second fundamental form ${\rm II}_{Z,z}$ is nondegenerate in the sense that its null space $${\rm Null}({\rm II}_{Z,z}) := \{ v \in T_z Z \mid \ {\rm II}_{Z,z}(v, w)
\mbox{ for all } w \in T_z Z \}$$ is zero. \end{itemize} Moreover, if  there exists a nonsingular point  $z \in Z$ satisfying the two equivalent conditions of (b),  \begin{itemize} \item[(c)] the third fundamental form
 ${\rm III}_{Z,z}$ of $Z$ at $z$ is, up to a nonzero scalar multiple, given by  a single cubic form  $f^z \in \Sym^3 T^*_z Z$, which satisfies $$f^z( \frac{\p}{\p y^i}, \frac{\p}{\p y^j}, \frac{\p}{\p y^k}) =  \frac{\partial^3 E(0)}{\partial y^i \partial y^j \partial y^k} $$ for all $1 \leq i, j, k \leq n$ in terms of the coordinates in Lemma \ref{l.coord};
 \item[(d)] the second fundamental form    ${\rm II}_{Z, z}$ is, up to a nonzero scalar multiple, given by the $n$-dimensional system of quadratic forms on $T_z Z$
     $$\{ f^z(v, \cdot, \cdot) \in \Sym^2 T^*_z Z \mid v \in T_z A\}$$ obtained by the contraction of the cubic form $f^z$ in (c); and  \item[(e)] $Z$ is linearly nondegenerate  in $\BP V$, namely, it is not contained in any hyperplane of $\BP V$.   \end{itemize} \end{proposition}

\begin{proof}
In terms of the local coordinates   in Lemma \ref{l.coord},
elements of $H^0(\BP V, \sO(1))$ whose restrictions to $Z$ vanish to the second order at $z$ are spanned by $F^k, 1 \leq k \leq n$ and $E$. Thus
 ${\rm II}_{Z,z}$ is given by  the $n+1$ quadratic forms
$$ \sum_{i,j=1}^n \frac{\partial^2 F^k(0)}{\partial y^i \partial y^j} {\rm d}y^i \cdot {\rm d} y^j \mbox{ for } 1 \leq k \leq n \mbox{ and }
\sum_{i,j=1}^n \frac{\partial^2 E(0)}{\partial y^i \partial y^j} {\rm d}y^i \cdot {\rm d} y^j.$$
But the latter vanishes by  (\ref{e.iv}).  By Lemma \ref{l.coord} (3), this implies (a).

To prove (b), note that ${\rm Im}({\rm II}_{Z,z}) = D_z/T_z Z$ if and only if the $n$ quadratic forms
$$ \sum_{i,j=1}^n \frac{\partial^2 F^k(0)}{\partial y^i \partial y^j} {\rm d}y^i \cdot {\rm d} y^j  \mbox{ for } 1 \leq k \leq n$$ are linearly independent. But  $$ \sum_{i,j, k=1}^n c^k \frac{\partial^2 F^k(0)}{\partial y^i \partial y^j} {\rm d}y^i \cdot {\rm d} y^j = 0 $$ for some $c^k \in \C, 1 \leq k \leq n,$ exactly when $\sum_{k=1}^n c^k \frac{\p}{\p y^k} \in T_z Z$ is in ${\rm Null}({\rm II}_{Z,z})$. This proves (b).

If  ${\rm Im}({\rm II}_{Z,z}) = D_z/T_z Z$, then
by Lemma \ref{l.coord}, the third fundamental form corresponds to the cubic form $f^z$ given by
$$ \sum_{i,j,k=1}^n \frac{\partial^3 E(0)}{\partial y^i \partial y^j \partial y^k} {\rm d}y^i \cdot {\rm d} y^j \cdot {\rm d} y^k.$$ This proves (c).

By  (\ref{e.i}) and (\ref{e.iii}), we have
$$ \sum_{i,j=1}^n \frac{\partial^2 F^k(0)}{\partial y^i \partial y^j} {\rm d}y^i \cdot {\rm d} y^j = - \sum_{i,j=1}^n \frac{\partial^3 E(0)}{\partial y^i \partial y^j \partial y^k} {\rm d}y^i \cdot {\rm d} y^j.$$ The left hand side is ${\rm II}_{Z,z}$ and the right hand side is the contraction of $-f^z$ by $\frac{\p}{\p y^k}$.  This proves (d).

Finally, in the setting of (c) and (d), the third fundamental form ${\rm III}_{z,z}$ is nonzero, which implies  (e).
\end{proof}

Now we look at nonsingular Legendrian varieties. We say that $Z \subset \BP V$ is nonlinear if it is not isomorphic to a linear subspace $\BP^{n} \subset \BP^{2n+1}$.

\begin{proposition}\label{p.LM}
Let $Z \subset \BP V$ be a  nonsingular and nonlinear  Legendrian variety. \begin{itemize} \item[(i)] A general point $z\in Z$ satisfies ${\rm Im}({\rm II}_{Z,z}) = D_z/T_z Z.$  \item[(ii)] For a general point $z \in Z$,
let
$${\rm Bs}({\rm II}_{Z,z}) := \{ v \in T_z Z \mid {\rm II}_{Z,z}(v,v) =0\}$$ be the base locus of ${\rm II}_{Z,z}$. Then it is exactly the set of tangent vectors to lines on $Z$ passing through $z$.
\item[(iii)]
The dual variety $Z^* \subset \BP V^*$ is a hypersurface in $\BP V^*$. \end{itemize} \end{proposition}

\begin{proof}
If $Z$ is nonsingular and not a linear subspace, then its Gauss map is birational (see Corollary on page 124 of \cite{FP}), which implies that ${\rm II}_{Z,z}$
is nondegenerate at a general point $z \in Z$ (see Proposition on page 111 of \cite{FP}).
Thus (i) is a consequence of Proposition \ref{p.III} (b).

(ii) is  Theorem 16 of \cite{LM}.

As noted in Proposition 17 of \cite{LM}, the dual variety $Z^*$ is isomorphic to the tangent variety ${\rm Tan}(Z) \subset \BP V$ via the isomorphism $V^* \cong V$ induced by the symplectic form $\sigma$. Thus to prove (iii),  it suffices to show that $\dim {\rm Tan}(Z)= 2 \cdot \dim Z.$ But if $\dim {\rm Tan}(Z) < 2 \cdot \dim Z,$ then ${\rm Tan}(Z) = {\rm Sec}(Z)$ by Corollary in page 123 of \cite{FP}, which implies that $Z$ has a degenerate secant variety. Then the third fundamental form of $Z$ at a general point vanishes identically by Theorem in page 135 of \cite{FP}, a contradiction to Proposition \ref{p.II} (c) and (d).
\end{proof}

\begin{definition}\label{d.Hessian} A cubic form $f \in \Sym^3 W^*$ on a vector space $W$ has {\em nonzero Hessian} if for some $w \in W,$ the quadratic form $f(w, \cdot, \cdot) \in \Sym^2 W^*$ is nondegenerate, namely, $$\{ u \in W \mid f(w, u, v) = 0 \mbox{ for all } v \in W\} =0.$$
\end{definition}

From Proposition \ref{p.LM}, we deduce the following.

\begin{proposition}\label{p.III}
Let $Z \subset \BP V$ and $z \in Z$ be as in Proposition \ref{p.LM}. Write $W = T_z Z$.   Then \begin{itemize} \item[(i)]
the cubic form $f^z$ in Proposition \ref{p.II} (c) has nonzero Hessian; and \item[(ii)] for the  cubic hypersurface $Y^z \subset \BP W$ determined by the cubic form $f^z$,
   its (set-theoretical) singular locus ${\rm Sing}(Y^z) \subset \BP W$ of $Y^z$ is nonsingular.
\end{itemize} \end{proposition}

\begin{proof}
By Proposition \ref{p.LM} (iii),  the dual variety of $Z$ is a hypersurface on $\BP V^*$, which implies that the system of quadratic forms ${\rm II}_{Z,z}$ at a general point $z \in Z$ contains a nondegenerate quadratic form by Proposition in page 112 of \cite{FP}. By
 Proposition \ref{p.II} (d), this implies that
 for a general $w \in T_z Z$, the contraction $f^z(w, \cdot, \cdot) \in \Sym^2 T^*_z Z$ is a nondegenerate quadratic form. Thus $f^z$ has nonzero Hessian, proving (i).

The affine cone in $W$ of ${\rm Sing}(Y^z) \subset \BP W$ is $$ \{ v \in W \mid f^z(v, v, u) = 0 \mbox{ for all } u \in W\}.$$ By Proposition \ref{p.II} (d), this is exactly the base locus ${\rm Bs}({\rm II}_{Z,z})$ in Proposition \ref{p.LM} (ii).
Recall the general fact that for a nonsingular projective variety $Z \subset \BP V$,
the variety $\sC_z \subset \BP T_z Z$ consisting of tangent directions to lines on $Z$ through $z$ is nonsingular for a general $z \in Z$ (e.g. by Proposition 1.5 of \cite{Hw01}).
In our setting, this $\sC_z$ is exactly  ${\rm Sing}(Y^z)$ by Proposition \ref{p.LM} (ii). This implies (ii).
\end{proof}

\section{Characterizing subadjoint varieties by their third fundamental forms}\label{s.subadjoint}

From Definition \ref{d.subadjoint},  the subadjoint variety associated to a simple Lie algebra $\fg$ is a homogeneous Legendrian variety $Z^{\fg} \subset \BP V^{\fg}$ in a symplectic vector space $V^{\fg}$ determined by $\fg$. It is well-known that their fundamental forms can be described as follows (see Corollary 26 of \cite{LM} and the references therein).

\begin{proposition}\label{p.EKP}
The third fundamental form of a  subadjoint variety $Z^{\fg} \subset \BP V^{\fg}$  is
isomorphic to  the determinant of a semisimple Jordan algebra of rank 3. For each $\fg$, the corresponding cubic hypersurface $Y^{\fg} \subset \BP W^{\fg},$ where $W^{\fg}$ denotes $T_z Z^{\fg}$ for a point $z \in Z^{\fg}$, can be described as follows. \begin{itemize}
\item[(i)] For $\fg= \fso_7$, the cubic form is isomorphic to $s^2 t$ in two variables
$(s,t)$. The singular locus of $Y^{\fg}$ is a single point in $\BP W^{\fg} \cong \BP^1.$
\item[(ii)] For $\fg = \fso_8$, the cubic hypersurface $Y^{\fg} \subset \BP^2$ is the union of three lines intersecting at three points in $\BP W^{\fg} \cong \BP^2$. Three intersections points,
    which are not collinear, are the singular locus of $Y^{\fg}$.
\item[(iii)] For $\fg= \fso_{n+5}, n \geq 4$, the cubic hypersurface $Y^{\fg}$ is the union of a hyperplane and a
quadric hypersurface with an isolated singularity outside the hyperplane.
\item[(iv)] For $\fg$ of type $G_2$, the cubic form is isomorphic to $s^3$ in one variable $(s)$. The hypersurface $Y^{\fg}$ is empty.
\item[(v)] For $\fg$ of type $F_4, E_6, E_7$ or $ E_8,$ the cubic hypersurface is the secant of the following  four Severi varieties:
     $$v_2(\BP^2) \subset \BP^5, \ \BP^2 \times \BP^2 \subset \BP^8, \
{\rm Gr}(2, 6) \subset \BP^{14}, \
  \mbox{ and } \
\BO \BP^2 \subset \BP^{26}. $$  \end{itemize} \end{proposition}

We reformulate Theorem \ref{t.III} as follows.

\begin{theorem}\label{t.HY}
Let $Z \subset \BP V$ be a nonsingular and nonlinear Legendrian variety of dimension $n \geq 2$.
Let $Y^z \subset \BP T_z Z$ be the cubic hypersurface defined by ${\rm III}_{Z,z}$ at a general point $z \in Z$. If $Y^z\subset \BP T_z Z$ is isomorphic to $Y^{\fg} \subset \BP W^{\fg}$ in Proposition \ref{p.EKP}, then $Z \subset \BP V$ is isomorphic to $Z^{\fg} \subset \BP V^{\fg}$. \end{theorem}

\begin{proof} Since $\dim Z = n \geq 2$, we can ignore  $\fg$ of type $G_2$, the case (iv) of Proposition \ref{p.EKP}. By the assumption, we have $\dim Z = \dim Z^{\fg}$ for some $\fg$ in Proposition \ref{p.EKP}

When $\fg $ is in the case (v) of Proposition \ref{p.EKP}, the list in Definition \ref{d.subadjoint} shows that the subadjoint variety is the minimal embedding of an irreducible Hermitian symmetric space, different from projective space or the hyperquadric. Since ${\rm II}_{Z,z}$ is determined by ${\rm III}_{Z,z}$, the projective variety $Z \subset \BP V$ has the same fundamental forms as those of the subadjoint variety $Z^{\fg}$. Thus $Z \subset \BP V$ is isomorphic to $Z^{\fg} \subset \BP V^{\fg}$ by the characterization of irreducible Hermitian symmetric spaces by their fundamental forms in Theorem 1 of \cite{HY}.

When $\fg$ is in the case (i) of Proposition \ref{p.EKP}, the Legendrian variety $Z$ is a surface and
the singular locus of the cubic defined by ${\rm III}_{Z,z}$ is a single point.
By  Proposition \ref{p.LM} (i), there exists a unique line through a general point of $Z$.
Since $Z$ is a nonsingular surface, this implies that $Z$ is a ruled surface. Thus $Z \cong \BP^1 \times \Q^1$ by Proposition 8 in \cite{LM}.

When $\fg$ is in the cases (ii) or (iii) of Proposition \ref{p.EKP},
we have lines covering $Z$ with $(n-3)$-dimensional nontrivial deformations fixing a general point $z \in Z$ corresponding to the $(n-3)$-dimensional component of the singular locus of $Y^{\fg}$ (two points in the case (ii) with $n=3$), which is isomorphic to  the quadric hypersurface $\Q^{n-3} \subset \BP^{n-2}$.  Projective varieties with such a large family of lines have been classified in Theorem 1.4 of \cite{LP}. Since the $(n-3)$-dimensional family of lines through a general point $z \in Z$ is isomorphic to the quadric hypersurface $\Q^{n-3} \subset \BP^{n-2},$ we see that $Z$ is the hyperquadric fibration (1.4.3.f) of \cite{LP}, namely, there is a surjective morphism $\psi: Z \to B$ to a curve $B$ whose general fiber is a quadric hypersurface. Furthermore, the singular locus of $Y^{\fg}$ has one additional isolated point. Thus we can use the next lemma to complete the proof. \end{proof}

\begin{lemma}\label{l.fusion}
Let $\Q^{n-1} \subset \BP^n$ be the nonsingular quadric hypersurface of dimension $n-1 \geq 1$.
Let $Z \subset \BP^{2n+1}$ be a linearly nondegenerate nonsingular projective variety of dimension $n$ equipped with a morphism $\psi: Z \to B$ to a projective curve $B$ satisfying the following properties. \begin{itemize} \item[(a)] For a general point $z \in Z$, there exists a linear subspace $\BP^n_z \subset \BP^{2n+1}$ such that the  fiber $Q_z := \psi^{-1}(\psi(z))$ of $\psi$ through $z$ is contained in  $\BP_z^n \cap Z$ and the inclusion $Q_z \subset \BP^n_z$ is isomorphic to $\Q^{n-1} \subset \BP^n$.  \item[(b)]   For a general point $z \in Z$, there exists a unique line $\ell_z \subset Z$ transversal to $Q_z$.  \end{itemize} Then $Z$ is isomorphic to the Segre embedding $$\BP^1 \times \Q^{n-1} \subset \BP (\C^2 \otimes \C^{n+1}) = \BP^{2n+1}.$$ \end{lemma}

\begin{proof} Because $\ell_z \cap Q_z = z$ by (a), the morphism $\psi$ sends $\ell_z$ isomorphically to $B$.
Let $\sL^o$ be the subset of the Grassmannian of lines on $\BP^{2n+1}$
parameterizing the lines $$\{\ell_z \subset Z \mid \mbox{ general } z \in Z \}$$
given in (b). Let $\sL$ be the closure of $\sL^o$ in the Grassmannian.
Then all members of $\sL$ are lines in $Z$ and must be sent  isomorphically to $B$ by $\psi$. So they  are transversal to fibers of $\psi$.
It follows that all fibers of $\psi$ must be nonsingular. A general $\ell_z$ in (b) must have trivial normal bundle in $Z$, otherwise there would be nontrivial deformations of $\ell_z$ in $Z$ fixing $z$, violating the uniqueness in  (b).  Thus we can find an algebraic subset $E \subset Z$ of codimension at least $2$, such that  members of $\sL$ determine a foliation of rank 1 on $Z\setminus E$,  which is transversal to the fibration $\psi$. It follows that this foliation is defined everywhere on $Z$ inducing a splitting of the tangent bundle of $TZ$ (see  Proposition 5 in \cite{Li}). This implies that $Z \cong \BP^1 \times \Q^{n-1}$. Since $Z$ is linearly nondegenerate, it must be the Segre embedding.
 \end{proof}

\begin{corollary}\label{c.III}
In Proposition \ref{p.III}, the cubic hypersurface $Y^z$ is irreducible and reduced, unless $Z$ is  isomorphic to the subadjoint variety $\BP^1 \times \Q^{n-1}$  associated with $\fg= \fso_{n+5}$. \end{corollary}

\begin{proof}
If  the cubic hypersurface is not reduced, the reduction (the underlying reduced hypersurface) must be a hyperplane. This is a contradiction to  Proposition \ref{p.III} (i). So the cubic hypersurface is reduced. Furthermore, if it is not irreducible, then it must be isomorphic to (ii) or (iii) in Proposition \ref{p.EKP}. Thus Corollary \ref{c.III} follows from Theorem \ref{t.HY}. \end{proof}

\section{Prolongations of infinitesimal automorphisms of cubic hypersurfaces}\label{s.cubic}

In this section, we recall  the main result of \cite{Hw20} with some supplement to fix a minor  gap in the argument.

\begin{definition}\label{d.prolong}
Let $S \subset \BP W$ be an irreducible variety and let  $\widehat{S} \subset W$ be the corresponding affine cone in the vector space $W$. For a nonsingular point $s \in \widehat{S}$, denote by $T_s \widehat{S} \subset W$ the affine tangent space at $s$.   \begin{itemize} \item[(i)] The Lie algebra $\aut(\widehat{S}) \subset {\rm End}(W)$ of infinitesimal automorphisms of $\widehat{S} \subset W$        consists of endomorphisms $ \varphi \in {\rm End}(W)$ satisfying $\varphi(s) \in T_s \widehat{S}$ for any nonsingular point $s \in \widehat{S}$. This is the Lie algebra of the linear automorphism group ${\rm Aut}(\widehat{S}) \subset {\rm GL}(W)$ of the affine cone, which is the inverse image of the linear automorphism group ${\rm Aut}(S) \subset {\rm PGL}(W)$ of $S \subset \BP V$  under the projection ${\rm GL}(W) \to {\rm PGL}(W)$.
\item[(ii)] For an element $A \in  \Hom (\Sym^2 W, W),$ denote by $A_{uv} = A_{vu} \in W$ its value
at $u, v \in W$. Then $A$ is a {\em prolongation} of $\aut(\widehat{S})$ if
for each $w \in W$, the endomorphism $A_w \in {\rm End}(W)$ defined by $A_w (u) := A_{wu}$ belongs to $\aut(\widehat{S})$. The vector space of all prolongations of $\aut(\widehat{S})$
 is denoted by $\aut(\widehat{S})^{(1)}.$ \end{itemize} \end{definition}

 The following is Theorem 2.1 of \cite{Hw20}, which is an easy consequence of the classification in \cite{FH18} of all linearly nondegenerate nonsingular subvariety $S \subset \BP W$ with $\aut(\widehat{S})^{(1)} \neq 0$.

\begin{theorem}\label{t.Severi}
 Let $S \subset \BP W$ be a linearly nondegenerate nonsingular   subvariety with $\aut(\widehat{S})^{(1)} \neq 0.$  If the secant variety ${\rm Sec}(S) \subset \BP W$ is a
hypersurface, then $S \subset \BP V$ is one of the following four Severi varieties:
 $$v_2(\BP^2) \subset \BP^5, \ \BP^2 \times \BP^2 \subset \BP^8, \
{\rm Gr}(2, 6) \subset \BP^{14}, \
  \mbox{ and } \
\BO \BP^2 \subset \BP^{26}. $$
 \end{theorem}

 The following  elementary fact is from Lemma 3.6 of \cite{Hw20}.

\begin{lemma}\label{l.3.6}
 Let
$Y \subset \BP W$ be an irreducible reduced cubic hypersurface defined by a cubic form $f \in \Sym^3 W^*$ on a vector space $W$. Then there  exists a linear functional  $\chi: \aut(\widehat{Y}) \to \C$ such that the Lie algebra $\aut(\widehat{Y}) \subset {\rm End}(W)$   consists of endomorphisms $ \varphi \in {\rm End}(W)$ satisfying $$ f(\varphi (u), v, w) + f(u, \varphi (v), w) + f(u, v, \varphi(w)) = \chi(\varphi) \ f(u, v, w) $$ for all $u,v,w \in W$. \end{lemma}

\begin{definition}\label{d.cubic} Let
$Y \subset \BP W$ be an irreducible reduced cubic hypersurface defined by a cubic form $f \in \Sym^3 W^*$ on a vector space $W$.
\begin{itemize}
\item[(i)] For $u, v \in W$, let $f_{uv} \in W^*$ be the linear functional defined by $f_{uv}(w) = f(u, v, w)$ for all $w \in W$.
    \item[(ii)] For $A \in \aut(\widehat{Y})^{(1)}$, let $\chi^A \in W^*$ be the linear functional defined by $\chi^A(u) = \chi(A_u)$ where $\chi: \aut(\widehat{Y}) \to \C$ is from Lemma \ref{l.3.6}.
\item[(iii)] For a complex number $a \in \C$, let  $\Xi^a_Y \subset \aut(\widehat{Y})^{(1)}$ be the subspace consisting of  $ A \in \aut(\widehat{Y})^{(1)}$ that satisfies $$ A_{uv}
=  a \chi^A(u) v + a \chi^A(v) u + h^A(f_{uv}) $$  for some $ h^A \in \Hom(W^*, W)$ and all $u, v \in W$.
 \end{itemize} \end{definition}

The following is the main result in this section and is exactly Theorem 1.6 of \cite{Hw20}.

\begin{theorem}\label{t.cubic} Let $Y \subset \BP W$ be an irreducible cubic hypersurface  with nonzero Hessian. Let ${\rm Sing}(Y) \subset Y$ be the singular locus of $Y$.
 Assume that \begin{itemize} \item[(a)] the  singular locus
${\rm Sing}(Y)$ is nonsingular; and \item[(b)]
$\Xi_Y^a \neq 0$ for some $a \neq \frac{1}{4}.$ \end{itemize}
Then ${\rm Sing}(Y)$ is one of the four Severi varieties (in Theorem \ref{t.Severi}) and  $Y$ is the secant variety ${\rm Sec}({\rm Sing}(Y))$ of the Severi variety.
\end{theorem}

There is a minor gap in the proof of Theorem \ref{t.cubic} in \cite{Hw20}. Let us explain this.
The key ingredients of the proof of Theorem \ref{t.cubic} are Theorem \ref{t.Severi} and   the following, which is Theorem 5.1 of \cite{Hw20}.

\begin{theorem}\label{t.5.1}
Let $Y \subset \BP W$ be an irreducible cubic hypersurface defined by a cubic form $f \in \Sym^3 W^*$ with nonzero Hessian.
 Assume that $\Xi_Y^a \neq 0$ for some $a \neq \frac{1}{4}.$  Then $Y = {\rm Sec}({\rm Sing}(Y)).$ \end{theorem}

Yewon Jeong has pointed out that the proof of  Theorem 5.1 of \cite{Hw20} (that is, Theorem \ref{t.5.1} above) requires $a \neq 0$ at the last line in page 41 of \cite{Hw20}. To fill this gap, we give  a proof of Theorem \ref{t.5.1} here, when $a=0$.  We use the following two lemmata.

The following lemma is from Lemma 3.4 (4) and Definition 3.5 of \cite{Hw20}.

\begin{lemma}\label{l.3.4}
Suppose that $f \in \Sym^3 W^*$ has nonzero Hessian. Then \begin{itemize} \item[(i)] the homomorphism $f_u: W \to W^*$ defined by $f_u(w) = f_{uw}$ is an isomorphism for a general $u \in W$; and \item[(ii)]   for any dense open subset $U \subset W$,  $$\bigcap_{ u \in U} {\rm Ker}(f_{uu})  = 0,$$ where ${\rm Ker}(f_{uu}) \subset W$ is the hyperplane given by $f_{uu} \in W^*$.  \end{itemize} \end{lemma}

In the next lemma,  (i) is a reformulation of Proposition 5.6 of \cite{Hw20} and (ii) is (C1) in page 41 of \cite{Hw20}.

\begin{lemma}\label{l.5.5}
For an irreducible cubic hypersurface $Y \subset \BP W$ with nonzero Hessian, if $Y \neq {\rm Sec}( {\rm Sing}(Y))$ and  $\Xi_Y^a \neq 0$ for some $a \neq \frac{1}{4}$, then there  exists a nonzero $A \in \Xi_Y^a$ such that \begin{itemize} \item[(i)] the one-dimensional subspace
$\C \cdot A \subset \Xi_Y^a$  is invariant under the connected Lie group ${\rm Aut}^o(\widehat{Y}) \subset {\rm GL}(W)$ with Lie algebra $\aut(\widehat{Y}) \subset \fgl(W) = {\rm End}(W)$; and \item[(ii)] for a general point $w \in \widehat{Y}$, the element $h^A(f_{ww})$ is contained in  $\Gamma_w \setminus B_w$, where $$\Gamma_w = \{ v \in W \mid f_{uw} \in \C \cdot f_{ww} \subset W^* \}$$
 is the Gauss fiber of $\widehat{Y}$ through $w$ (by Proposition 4.2 of \cite{Hw20}) and   $B_w = \Gamma_w \cap {\rm Ker}(\chi^A)$ (from Proposition 5.3 of \cite{Hw20}). \end{itemize} \end{lemma}

\begin{proof}[Proof of Theorem \ref{t.5.1} when $a=0$.]
We may assume the setting of  Lemma \ref{l.5.5} and derive a contradiction. Let  $A \in \Xi^a_Y$ with $a=0$ be as in Lemma \ref{l.5.5}. Then
\begin{equation}\label{e.def}
A_{uv} = h^A ( f_{uv}) \ \mbox{ for all }u, v \in W.
\end{equation} Lemma \ref{l.5.5} (ii) implies
 \begin{equation}\label{e.imker} {\rm Im}(h^A) \not\subset {\rm Ker}(\chi^A). \end{equation}
Since $f_u: W \to W^*$  is an isomorphism for a general $u \in W$ by  Lemma \ref{l.3.4}, (\ref{e.def}) implies \begin{equation}\label{e.im} {\rm Im}(A_u) = {\rm Im}(h^A) \ \mbox{ for all general } u \in W. \end{equation}
Lemma \ref{l.3.6} with $\varphi = A_w$
 gives $$ 3 f_{uu} (A_u(w)) = 3 f(A_{uw}, u, u) = 3 f(A_w(u), u, u) = \chi^A(w) \ f(u, u, u)$$ for any $u, w \in W$.  This shows  $A_u({\rm Ker}(\chi^A)) \subset {\rm Ker}(f_{uu}).$ Moreover, if $u \not\in \widehat{Y}$, then ${\rm Ker}(\chi^A) = A_u^{-1}({\rm Ker}(f_{uu})). $
Together with (\ref{e.im}), we obtain \begin{equation}\label{e.ker} A_u({\rm Ker}(\chi^A)) = {\rm Im}(A_u) \cap {\rm Ker}(f_{uu}) =  {\rm Im}(h^A) \cap {\rm Ker}(f_{uu})\end{equation} for all general $u \in W$.
Note that $$A_u({\rm Ker} (\chi^A)) \ \subset \ {\rm Ker}(\chi^A)$$
 by the condition (i) of Lemma \ref{l.5.5} and  $A_u \in \aut(\widehat{Y})$.
Thus \begin{equation}\label{e.subset} A_u({\rm Ker}(\chi^A)) \subset {\rm Im}(h^A) \cap {\rm Ker}(\chi^A).\end{equation}
By (\ref{e.imker}), the righthand side of (\ref{e.subset}) is a hyperplane in ${\rm Im}(h^A)$ containing (\ref{e.ker}). Consequently, $${\rm Im}(h^A) \cap {\rm Ker}(\chi^A) = {\rm Im}(h^A) \cap {\rm Ker}(f_{uu})$$ for all general $u \in W$. It follows that ${\rm Im}(h^A) \cap {\rm Ker}(\chi^A)$ is contained in the intersection of ${\rm Ker}(f_{uu})$ for all general $u \in W$, which is zero by Lemma \ref{l.3.4}. Consequently, the dimension of ${\rm Im}(h^A)$ is at most one. Then Lemma \ref{l.5.5} (ii) says that $\widehat{Y}$ has only one Gauss fiber, a contradiction.
\end{proof}

\section{Jets of contact vector fields tangent to  Legendrian varieties}\label{s.jet}

In this section, we will study  contact vector fields on $\BP V$,
namely, vector fields which generate 1-parameter families of automorphisms of $\BP V$ preserving
the contact structure $D \subset T \BP V$. We recall the  following general result on contact vector fields on contact manifolds, from Theorem 7.1 in Chapter I of \cite{Ko}.

 \begin{theorem}\label{t.Kobayashi}
 Let $M$ be a complex manifold with a contact structure $D \subset TM$ and the contact line bundle $L = TM/D$.
 Denote by $\aut(M, D)$ the Lie algebra of contact vector fields on $M$.
 Then the quotient homomorphism $\aut(M, D) \subset H^0(M, TM) \to H^0(M, L),$
 induces  a linear isomorphism  $\phi: \aut(M, D) \cong  H^0(M, L).$
 \end{theorem}

\begin{proposition}\label{p.Qformula}
Let $(V, \sigma)$ be a symplectic vector space with $\dim V = 2n +2$ and let $D \subset T \BP V$ be the contact structure determined by $\sigma$ with the line bundle $L = T \BP V/ D \cong \sO(2)$ from Definition \ref{d.Legendrian}.
The isomorphism $\phi$ in Theorem \ref{t.Kobayashi} between  the space of contact vector fields
on $\BP V$ and $\Sym^2 V^* = H^0(\BP V, L)$ can be explicitly given in linear coordinates as follows.
  Choose an affine cell $O\subset \BP V$
equipped with an inhomogeneous coordinate system $(x^1, \ldots, x^n, x^{n+1}, \ldots, x^{2n}, x^{2n+1})$ such that
a contact form is given by $$\theta = \sum_{k=1}^n (x^{n+k} {\rm d} x^k - x^k {\rm d} x^{n+k}) - {\rm d} x^{2n+1}$$ as in Lemma \ref{l.coord}.
 An element $Q \in \Sym^2 V^*$ can be written as a polynomial $q(x^1, \ldots, x^{2n+1})$ of degree at most $2.$ Then the corresponding contact vector field $\vec{Q}$ under the isomorphism $\phi$ of  Theorem \ref{t.Kobayashi} can be written on $O$ as
 \begin{eqnarray*} \vec{Q}|_O &=& \frac{1}{2} \sum_{k=1}^n(\frac{\partial q}{\partial x^{n+k}} - \frac{\partial q}{\partial x^{2n+1}} x^k) \frac{\partial}{\partial x^k}\\ & & +\frac{1}{2} \sum_{k=1}^n (- \frac{\partial q}{\partial x^k} - \frac{\partial q}{\partial x^{2n+1}} x^{n+k}) \frac{\partial}{\partial x^{n+k}}\\ & &  + \left( \frac{1}{2} \sum_{k=1}^n (\frac{\partial q}{\partial x^k} x^k + \frac{\partial q}{\partial x^{n+k}} x^{n+k}) - q \right) \frac{\partial}{\partial x^{2n+1}}. \end{eqnarray*} \end{proposition}

 \begin{proof}
 It is straightforward to check (see also Theorem 2.2 in \cite{Hw19}) that the above expression $\vec{Q}|_O$  satisfies $$\theta( \vec{Q}|_O) = q \mbox{ and } \ {\rm Lie}_{\vec{Q}|_O} \theta = - \frac{\partial q}{\partial x^{2n+1}} \theta.$$ So it is a contact vector field. It remains to check that this expression $\vec{Q}|_O$ gives a holomorphic vector field on $\BP V$. Using Euler's formula $$\sum_{k=1}^{2n+1} x^k \frac{\partial }{\partial x^k} h(x^1, \ldots, x^{2n+1}) = r \ h(x^1, \ldots, x^{2n+1})$$ for any homogeneous polynomial $h$ of degree $r$, we obtain  $$\frac{1}{2} \sum_{k=1}^n (\frac{\partial q}{\partial x^k} x^k + \frac{\partial q}{\partial x^{n+k}} x^{n+k}) -q = - \frac{1}{2} \frac{\partial q}{\partial x^{2n+1}} x^{2n+1} + b(x)$$ for some  polynomial $b(x) = b(x^1, \ldots, x^{2n+1}) $ of degree at most $ 1$. Thus the expression of $\vec{Q}|_O$ is reduced to
  \begin{equation}\label{e.b} b_0(x)  \sum_{k=1}^{2n+1} x^k \frac{\partial}{\partial x^k} + \sum_{k=1}^{2n+1} b_k(x)  \frac{\partial}{\partial x^k} \end{equation} for some polynomials $b_k(x), 0 \leq k \leq 2n+1,$ of degree at most $1$. It is easy to check that a vector field of  the form (\ref{e.b}) on $O$ can be extended to  a holomorphic vector field on $\BP V$. In fact, choose another inhomogeneous coordinate system $(y^1, \ldots, y^{2n+1})$ on an affine cell $O' \subset \BP V$ satisfying
 $$x^1 = \frac{1}{y^1}, \ x^2 = \frac{y^2}{y^1},  \ldots, \ x^{2n+1} = \frac{y^{2n+1}}{ y^1}.$$
 From $$ \frac{\partial}{\partial x^1} = - y^1(y^1 \frac{\partial}{\partial y^1} + \cdots + y^{2n+1} \frac{\partial}{\partial y^{2n+1}}) \mbox{ and } \frac{\partial}{\partial x^i} = y^1 \frac{\partial}{\partial y^i} \mbox{ for } 2 \leq i \leq 2n+1,$$  we see that  (\ref{e.b}) remains holomorphic in the coordinate system $(y^1, \ldots, y^{2n+1})$. Hence, it is holomorphic on $O \cup O'$. Since the complement $\BP V \setminus (O \cup O')$ has codimension 2 in $\BP V$, we can extend (\ref{e.b}) to a holomorphic vector field on $\BP V$. \end{proof}

\begin{proposition}\label{p.jet}
Let $Z \subset \BP V$ be a Legendrian variety and
let $z \in Z$ be a nonsingular  point satisfying the condition in  Proposition \ref{p.II} (b).  Write $W = T_{z}Z$ and let $f^z \in \Sym^3 W^*$ be the nonzero cubic form in Proposition \ref{p.II} (c).
Let $\vec{Q}$  be a contact vector field on $\BP V$ that is tangent to $Z$ and write $\vec{Q}|_Z$ for the vector field on $Z$ given by the restriction. Assume that \begin{itemize} \item[(a)]   the cubic hypersurface $Y^z \subset \BP W$ defined by $f^z$ is irreducible and reduced; and \item[(b)] the vector field  $\vec{Q}|_Z$ vanishes to the second order at $z$, namely, its linear part at $z$ vanishes.  \end{itemize} Then the 2-jet of $\vec{Q}|_Z$ at $z$ determines an element
$A \in \Hom(\Sym^2 W, W)$, which belongs to $\aut(\widehat{Y}^z)^{(1)}.$ \end{proposition}

\begin{proof}
The 1-parameter family of automorphisms of $Z$ generated by $\vec{Q}$ preserve the closure of the subset $$\bigcup_{\mbox{ general } z \in Z} Y^z  \ \subset \ \BP TZ.$$
Thus the 2-jet of $\vec{Q}|_Z$ belongs to $\aut(\widehat{Y}^z)^{(1)}$ by Proposition 1.2.1 of \cite{HM} (see also Proposition 5.9 of \cite{FH12}). \end{proof}

The main result of this section is the following.

\begin{theorem}\label{t.vecQ} Assume the setting of Proposition  \ref{p.jet}.
\begin{itemize} \item[(i)]
The 2-jet $A \in \aut(\widehat{Y})^{(1)}$ of $\vec{Q}|_Z$  at $z$ belongs to $\Xi^{\frac{1}{2}}_Y, $ namely,
$$ A_{uv} = \frac{1}{2}\chi^A (u) v + \frac{1}{2} \chi^A (v) u + h (f^z_{uv}) \mbox{ for all } u, v \in W$$
for some  $h \in \Hom(W^*, W)$.  \item[(ii)] Assume furthermore that $f^z$ has nonzero Hessian. If $A=0$, namely, if  the vector field $\vec{Q}|_Z$ vanishes to the third order at $z$, then $\vec{Q} =0$ on $\BP V$. \end{itemize} \end{theorem}

 \begin{proof}
 Let us use the  coordinates from Lemma  \ref{l.coord}.
   Let $q(x^1, \ldots, x^{2n+1})$ be the polynomial of degree at most $2$ in Proposition \ref{p.Qformula}, corresponding to the vector field $\vec{Q}$. From the expression of $\vec{Q}|_O$ given in Proposition \ref{p.Qformula}, the vanishing of $\vec{Q}$ at $z$ implies
\begin{equation}\label{e.vanish1}
q(0) = \frac{\p q (0)}{\p x^k} = \frac{\p q(0)}{\p x^{n+k}} = 0 \end{equation} for all $1 \leq k \leq n.$

   Let us write $$\vec{Q}|_{U\cap Z} = \sum_{k=1}^n B^k(y) \frac{\partial}{\partial y^k}$$ for some holomorphic functions $B^k(y^1, \ldots, y^n)$ defined on $U\cap Z$.
Then $$B^k = \vec{Q}|_{U\cap Z} (y^k) =  \vec{Q}(x^k)|_{U \cap Z} = \frac{1}{2} (\frac{\partial q}{\partial x^{n+k}} - \frac{\partial q}{\partial x^{2n+1}} x^k) |_{U \cap Z}.$$ Applying chain rule and (\ref{e.i}), we obtain
\begin{eqnarray*} 2 \frac{\partial B^k}{\partial y^m} &=& \sum_{j=1}^n \frac{\partial}{\partial x^j}\left( \frac{\partial q}{\partial x^{n+k}} - \frac{\partial q}{\partial x^{2n+1}} x^k \right) \frac{\partial x^j}{\partial y^m} \\ & & + \sum_{j=1}^n \frac{\partial}{\partial x^{n+j}} \left( \frac{\partial q}{\partial x^{n+k}} - \frac{\partial q}{\partial x^{2n+1}} x^k \right) \frac{\partial x^{n+j}}{\partial y^m} \\ && + \frac{\partial}{\partial x^{2n+1}} \left( \frac{\partial q}{\partial x^{n+k}} - \frac{\partial q}{\partial x^{2n+1}} x^k \right) \frac{\partial x^{2n+1}}{\partial y^m} \\ &=&
\frac{\partial^2 q}{\partial x^m \partial x^{n+k}} - \frac{\partial^2 q}{\partial x^m \partial x^{2n+1}} x^k - \frac{\partial q}{\partial x^{2n+1}} \delta_{km} \\ & & +
\sum_{j=1}^n \left( \frac{\partial^2 q}{\partial x^{n+j} \partial x^{n+k} }-\frac{ \partial^2 q}{\partial x^{n+j} \partial x^{2n+1}} x^k \right) \frac{\partial^2 F}{\partial y^m \partial y^j} \\ & & + \left( \frac{\partial^2 q}{\partial x^{2n+1} \partial x^{n+k}} - \frac{\partial^2 q}{\partial (x^{n+1})^2} x^k \right) \frac{\partial E}{\partial y^m}. \end{eqnarray*}
By our assumption that the vector field $\vec{Q}|_Z$ vanishes to
the  second order at $z$, $$ \frac{\p B^k (0)}{\p y^m} =0 \mbox{  for all } 1 \leq k, m \leq n.$$ By (\ref{e.ii}), this implies \begin{equation}\label{e.vanish2}
\frac{\p^2 q(0)}{\p x^m \p x^{n+k}} = \frac{\p q (0)}{\p x^{2n+1} } \delta_{km} \mbox{ for all } 1 \leq k, m \leq n. \end{equation}

Taking derivative one more time, we have
\begin{eqnarray*} 2 \frac{\partial^2 B^k}{\partial y^{\ell} \partial y^m} & = & \frac{\partial}{\partial y^{\ell}}
\left( \frac{\partial^2 q}{\partial x^m \partial x^{n+k}} - \frac{\partial^2 q}{\partial x^m \partial x^{2n+1}} x^k - \frac{\partial q}{\partial x^{2n+1} } \delta_{km} \right) \\ & & + \sum_{j=1}^n  \frac{\partial}{\partial y^{\ell}} \left( \frac{\partial^2 q}{\partial x^{n+j} \partial x^{n+k} }-\frac{ \partial^2 q}{\partial x^{n+j} \partial x^{2n+1}} x^k \right) \cdot  \frac{\partial^2 F}{\partial y^m \partial y^j} \\ & & + \sum_{j=1}^n  \left( \frac{\partial^2 q}{\partial x^{n+j} \partial x^{n+k} }-\frac{ \partial^2 q}{\partial x^{n+j} \partial x^{2n+1}} x^k \right)   \frac{\partial^3 F}{ \partial y^{\ell} \partial y^m \partial y^j} \\ & & +  \frac{\partial}{\partial y^{\ell} } \left( \frac{\partial^2 q}{\partial x^{2n+1} \partial x^{n+k}} - \frac{\partial^2 q}{\partial (x^{n+1})^2} x^k \right) \cdot  \frac{\partial E}{\partial y^m} \\ & & + \left( \frac{\partial^2 q}{\partial x^{2n+1} \partial x^{n+k}} - \frac{\partial^2 q}{\partial (x^{n+1})^2} x^k \right) \frac{\partial^2 E}{\partial y^{\ell} \partial y^m}. \end{eqnarray*}
When we evaluate this at the point $z= (y^1= \ldots = y^n =0)$, the second, the fourth and the fifth lines of the righthand side vanish by (\ref{e.ii}) and (\ref{e.iv}). Since $q$ is a polynomial of degree at most $2$ in $x^1, \ldots, x^{2n+1}$, its third-order derivatives in $x^1, \ldots, x^{2n+1}$ must vanish identically. Thus we can write, using (\ref{e.iii}),
\begin{eqnarray}\label{e.B2}
  \frac{\partial^2 B^k(0)}{\partial y^{\ell} \partial y^m} & = & -\frac{1}{2} \frac{\partial^2 q (0) }{\partial x^{2n+1} \partial x^m} \delta_{k\ell} - \frac{1}{2} \frac{\partial^2 q (0)}{\partial x^{2n+1} \partial x^{\ell}} \delta_{km} \\ \nonumber &  &  - \frac{1}{2} \sum_{j=1}^n \frac{\partial^2 q (0) }{\partial x^{n+j} \partial x^{n+k}} \frac{\partial^3 E (0) }{\partial y^j \partial y^{\ell} \partial y^m}. \end{eqnarray}
  Define $\nu \in W^*$ by \begin{equation}\label{e.nu} \nu(\frac{\p}{\p y^i}) : = - \frac{1}{2} \frac{\p^2 q (0)}{\p x^{2n+1} \p x^i}\end{equation} and $h \in \Hom(W^*, W)$ by
\begin{equation}\label{e.h} h( {\rm d} y^i) := - \frac{1}{2} \sum_{k=1}^n \frac{\p^2 q (0)}{\p x^{n+i} \p x^{n+k}} \frac{\p}{\p y^k}.\end{equation} Then by  Proposition \ref{p.II} (c) and (\ref{e.B2}),
the 2-jet $A \in \Hom(\Sym^2 W, W)$ of $\vec{Q}|_Z$ at $z$ determined by
 $$ A_{\frac{\p}{\p y^{\ell}} \frac{\p}{\p y^m}}   :=  \sum_{k=1}^n  \frac{\partial^2 B^k (0)}{\partial y^{\ell} \partial y^m}  \frac{\p}{\p y_k}$$ satisfies
 \begin{equation}\label{e.star} A_{\frac{\p}{\p y^{\ell}} \frac{\p}{\p y^m}} = \nu(\frac{\p}{\p y^{\ell}})\frac{\p}{\p y^m} + \nu(v) \frac{\p}{\p y^{\ell}} + h(f^z_{\frac{\p}{\p y^{\ell}}\frac{\p}{\p y^m}}) \end{equation} for all $ 1 \leq \ell, m \leq n.$ To prove (i), it remains to show that \begin{equation}\label{e.nuchi} \nu = \frac{1}{2} \chi^A. \end{equation}

By (\ref{e.vanish1}) and (\ref{e.vanish2}),  the polynomial $q$ must be of the form
\begin{eqnarray*}
q(x^1, \ldots, x^{2n+1}) & =&  a x^{2n+1} + \sum_{i,j=1}^n b_{ij} x^i x^j + \sum_{i,j=1}^n c_{ij} x^{n+i} x^{n+j}\\ & &  + \sum_{i=1}^n d_i x^i x^{2n+1} + \sum_{i=1}^n e_i x^{n+i} x^{2n+1} \\ & &  + a \sum_{i=1}^n x^i x^{n+i} + g (x^{n+1})^2 \end{eqnarray*} for some complex numbers $a, b_{ij}= b_{ji}, c_{ij}=c_{ji}, d_{i}, e_i$ and $g$.  Then  (\ref{e.nu}) and (\ref{e.h}) become
  \begin{equation}\label{e.nuh} \nu(\frac{\p}{\p y^i}) = - \frac{1}{2} d_i \ \mbox{ and } \  h({\rm d} y^i) = - \sum_{k=1}^n c_{ik} \frac{\p}{\p y^k}. \end{equation}

 Let us write $ A_{\frac{\p}{\p y^{\ell}}\frac{\p}{\p y^m}} = \sum_{k=1}^n A^k_{\ell m} \frac{\p}{\p y^k}$ and $f^z_{\frac{\p}{\p y^i} \frac{\p}{\p y^j} \frac{\p}{\p y^k}} = f_{ijk}$. Then (\ref{e.star}) and  (\ref{e.nuh}) give $$  A^k_{\ell m}  = - \frac{1}{2}d_m \delta^k_{\ell} - \frac{1}{2} d_{\ell} \delta^k_m - \sum_{j=1}^n c_{jk} f_{\ell m j}.$$ Consequently,  $$ A^k_{ii}  = - d_i \delta^k_{i}  - \sum_{j=1}^n c_{jk} f_{j ii}$$ for all $1 \leq i, k \leq n$, which gives  \begin{equation}\label{e.Akii}
  \sum_{k=1}^n A^k_{ii}f_{kii}  = - d_i f_{iii}  - \sum_{j,k=1}^n c_{jk} f_{kii} f_{j ii}. \end{equation}
Lemma \ref{l.3.6} and
 $A \in \aut(\widehat{Y})^{(1)}$ imply
$$\sum_{\ell =1}^n \left( A^{\ell}_{mi} f_{\ell jk} + A^{\ell}_{mj} f_{i\ell k} + A^{\ell}_{mk}  f_{ij \ell}\right)  = \chi^A(\frac{\partial}{\partial y^m}) f_{ijk}.$$
Setting $i=j=k=m,$ we obtain
$$3 \sum_{\ell =1}^n A^{\ell}_{ii} f_{\ell ii} = \chi^A(\frac{\p}{\p y^i}) f_{iii}. $$
By (\ref{e.Akii}), this equation becomes
 \begin{equation}\label{e.Fiii}\chi^A(\frac{\partial}{\partial y^i}) f_{iii}  =   -3 d_i f_{iii} - 3 \sum_{j, k=1}^n  c_{jk}  f_{kii} f_{jii}\end{equation} for any $1 \leq i \leq n$.

Let $\sQ(y)$ be the holomorphic function $q|_{U \cap Z}$. By
putting the equations of $Z$ in Lemma \ref{l.coord} into $q(x^1, \ldots, x^{2n+1}),$ we obtain
\begin{eqnarray*}
 \sQ(y^1, \ldots, y^n) & := & a E + \sum_{i,j=1}^n b_{ij} y^i y^j + \sum_{i,j=1}^n c_{ij} \frac{\partial F}{\partial y^i} \frac{\partial F}{\partial y^j} \\ & & + \sum_{i=1}^n d_i y^i E + \sum_{i=1}^n e_i \frac{\partial F}{\partial y^i} E + a \sum_{i=1}^n y^i \frac{\partial F}{\partial y^i} + g E^2. \end{eqnarray*}
 From the relation of $Q$ and $\vec{Q}$ in Theorem \ref{t.Kobayashi} and Proposition \ref{p.Qformula}, the inclusion $TZ \subset D|_Z$ implies that  $Z$ is contained in the zero set of $q$. Thus the holomorphic function $\sQ(y)$ must be identically zero.
From (\ref{e.ii}), (\ref{e.iv}) and $\frac{\partial^2 \sQ (0)}{\partial y^i \partial y^j}= 0$, we obtain $b_{ij} = 0$ for all $1 \leq i, j \leq n$.
Using (\ref{e.ii}), (\ref{e.iv}) and
$$\frac{\partial^3}{\partial y^k \partial y^{\ell} \partial y^m}\left( \sum_{i=1}^n y^i \frac{\partial F}{\partial y^i} \right) = 3 \frac{\partial^3 F}{\partial y^k \partial y^{\ell} \partial y^m} + \sum_{i=1}^n y^i \frac{\partial^4 F}{\partial y^i \partial y^k \partial y^{\ell} \partial y^m},$$
we have
$$0 = \frac{\partial^3 \sQ (0)}{\partial y^k \partial y^{\ell} \partial y^m} = a \frac{\partial^3 E(0)}{\partial y^k \partial y^{\ell} \partial y^m} + 3 a \frac{\partial^3 F(0)}{\partial y^k \partial y^{\ell} \partial y^m}$$ for all $1\leq k, \ell, m \leq n$. Combining this  with (\ref{e.iii}) and $f^z \neq 0$, we obtain $a=0$.
So we are left with \begin{eqnarray}\label{e.sQ}
 \sQ(y) & =& \sum_{i,j=1}^n c_{ij} \frac{\partial F}{\partial y^i} \frac{\partial F}{\partial y^j} + \sum_{i=1}^n d_i y^i E + \sum_{i=1}^n e_i \frac{\partial F}{\partial y^i} E + g E^2.  \end{eqnarray}
By (\ref{e.ii}) and (\ref{e.iv}), if  we take the fourth derivative with respect to $y^k$ of the left hand side of (\ref{e.sQ}) and evaluate it at $y=0$, we are left with
$$
\frac{\p^4 \sQ(0)}{\p(y^k)^4} =
6 \sum_{i,j=1}^n c_{ij} \frac{\partial^3 E(0)}{\partial y^i\partial y^k \partial y^k } \frac{\partial^3 E(0)}{\partial y^j \partial y^k \partial y^k } + 4 d_k  \frac{\partial^3 E(0)}{\partial y^k\partial y^k \partial y^k } =0. $$
In other words, $$3 \sum_{j,k=1}^n c_{jk} f_{kii}f_{jii} + 2 d_i f_{iii} =0$$
for all $1 \leq i \leq n$. Combining it with (\ref{e.Fiii}), we have $\chi^A(\frac{\p}{\p y^i}) = - d_i = 2 \nu(\frac{\p}{\p y^i})$. This verifies (\ref{e.nuchi}), proving (i).

For (ii), recall that
Proposition 3.7 of \cite{Hw20} says that the association $A \mapsto \chi^A$ is injective if $f$ has nonzero Hessian. Thus $A=0$   implies $\chi^A=0 = h$. By (\ref{e.nuh}),   we have $d_i =0$ and $c_{ij} =0$ for all $1 \leq i, j \leq n$.   Thus (\ref{e.sQ}) yields $$ \sQ(y) = \sum_{k=1}^n e_k \frac{\partial F}{\partial x^k} E + g E^2 = 0.$$ Since $E(y)$ is not identically zero, we must have $\sum_{k=1}^n e_k \frac{\partial F}{\partial x^k} + g E =0.$ This means that $Z$ satisfies the linear
equation $\sum_{k=1} e_k x^{n+k} + g x^{2n+1} = 0$. Since $Z$ is linearly nondegenerate by Proposition \ref{p.II} (e), we must have $e_k = g=0$. Thus $q=0$,  proving (ii).
 \end{proof}

The following is Theorem 2 of \cite{Bu07} (see also Proposition 3 of \cite{Hw22} for a simple proof).

\begin{proposition}\label{p.Bu07}
Let $Z \subset \BP V$ be a nonsingular and nonlinear Legendrian variety. Then the image of $\aut(\widehat{Z}) \subset \fgl(V)$ in $H^0(\BP V, T \BP V)$ is contained in $\aut(\BP V, D) \cong \Sym^2 V^*$. \end{proposition}

\begin{proof}[Proof of  Theorem \ref{t.main}]
Assume that $\dim {\rm Ker}(\iota_z) \neq 0.$ By Proposition \ref{p.Bu07},
we have a contact vector field $\vec{Q}$ tangent to $Z$ vanishing to the second order at $z$.
If the cubic hypersurface $Y^z \subset \BP T_z(Z)$ defined by the third fundamental form is reducible, we know that $Z \cong \BP^1 \times \Q^{n-1}$ by Corollary \ref{c.III}. Thus we may assume that $Y^z$ is irreducible.
Since $f^z$ has nonzero Hessian by Proposition \ref{p.III}, Theorem \ref{t.vecQ} gives $\Xi^{\frac{1}{2}}_{Y^z} \neq 0$. As ${\rm Sing}(Y^z)$ is nonsingular by Proposition \ref{p.III},  we see that $Y^z$ is the secant variety of
a Severi variety by Theorem \ref{t.cubic}. It follows that $Z$ is a subadjoint variety by Theorem \ref{t.HY}.

Conversely, if $Z$ is a subadjoint variety, the embedding $Z \subset \BP V$ is equivariant, namely, all vector fields on $Z$ come from $\aut(\widehat{Z})$. As $Z$ is a Hermitian symmetric space, it is well-known (e.g., Section 3 of \cite{HY}) that for any given $z \in Z$, there exists a nonzero vector field on $Z$ vanishing to the second order at  $z$.  It follows that $\dim {\rm Ker}(\iota_z) \neq 0$.  \end{proof}

\bigskip
Jun-Muk Hwang(jmhwang@ibs.re.kr)

\smallskip

Center for Complex Geometry,
Institute for Basic Science (IBS),
Daejeon 34126, Republic of Korea

\end{document}